\def\bbbq{{\mathchoice 
{\setbox0=\hbox {$\displaystyle\rm Q$}\hbox
{\raise0.15\ht0\hbox to0pt{\kern0.4\wd0\vrule height0.8\ht0\hss}\box0}}
{\setbox0=\hbox {$\textstyle\rm Q$}\hbox
{\raise0.15\ht0\hbox to0pt{\kern0.4\wd0\vrule height0.8\ht0\hss}\box0}}
{\setbox0=\hbox {$\scriptstyle\rm Q$}\hbox
{\raise0.15\ht0\hbox to0pt{\kern0.4\wd0\vrule height0.7\ht0\hss}\box0}}
{\setbox0=\hbox {$\scriptscriptstyle\rm Q$}\hbox
{\raise0.15\ht0\hbox to0pt{\kern0.4\wd0\vrule height0.7\ht0\hss}\box0}}
}}

\def\bbbc{{\mathchoice 
{\setbox0=\hbox {$\displaystyle\rm C$}\hbox
{\hbox to0pt{\kern0.4\wd0\vrule height0.9\ht0\hss}\box0}}
{\setbox0=\hbox {$\textstyle\rm C$}\hbox
{\hbox to0pt{\kern0.4\wd0\vrule height0.9\ht0\hss}\box0}}
{\setbox0=\hbox {$\scriptstyle\rm C$}\hbox
{\hbox to0pt{\kern0.4\wd0\vrule height0.9\ht0\hss}\box0}}
{\setbox0=\hbox {$\scriptscriptstyle\rm C$}\hbox
{\hbox to0pt{\kern0.4\wd0\vrule height0.9\ht0\hss}\box0}}
}}

\font\fivesans=cmss10 at 5pt 
\font\sevensans=cmss10 at 7pt 
\font\tensans=cmss10   
\newfam\sansfam  
\textfont\sansfam=\tensans\scriptfont\sansfam=\sevensans  
\scriptscriptfont\sansfam=\fivesans  
\def\sans{\fam\sansfam\tensans} 
\def\bbbz{{\mathchoice {\hbox{$\sans\textstyle Z\kern-0.4em Z$}}  
{\hbox{$\sans\textstyle Z\kern-0.4em Z$}}  
{\hbox{$\sans\scriptstyle Z\kern-0.3em Z$}}  
{\hbox{$\sans\scriptscriptstyle Z\kern-0.2em Z$}}}} 
   
\def\slash#1{#1\kern-0.65em /}
\def\dirac{{\raise0.09em\hbox{/}}\kern-0.58em\partial}
\def\Dirac{{\raise0.09em\hbox{/}}\kern-0.69em D}

\def\m@th{\mathsurround=0pt}
\def\eqalign#1{\null\,\vcenter{\openup 3pt \m@th
\ialign{\strut\hfil$\displaystyle{##}$&$\displaystyle{{}##}$\hfil
\crcr#1\crcr}}\,}

\documentstyle[12pt]{article}
\begin{document}

\title{Differential Calculi and Linear Connections}

\author{Aristophanes Dimakis\\
        Institut f\"ur Theoretische Physik\\
        Bunsenstra\ss e 9, D-37073 G\"ottingen
\and    J. Madore \\
        Laboratoire de Physique Th\'eorique et Hautes 
        Energies\thanks{Laboratoire associ\'e au CNRS, {\rm URA D0063}}\\
        Universit\'e de Paris-Sud, B\^at. 211, F-91405 Orsay
       }
\date{January, 1996} 
\maketitle

\abstract{
A method is proposed for defining an arbitrary number of differential
calculi over a given noncommutative associative algebra.  As an example
the generalized quantum plane is studied.  It is found that there is a
strong correlation, but not a one-to-one correspondence, between the
module structure of the 1-forms and the metric torsion-free connections
on it. In the commutative limit the connection remains as a shadow of
the algebraic structure of the 1-forms.
}

\vfill
\noindent
LPTHE Orsay 95/79\par
\medskip
\eject

\section{Introduction}

We propose a method of defining an arbitrary number of differential
calculi over a given noncommutative associative algebra.  We shall
especially be interested in the generalized quantum plane, an algebra
which has a commutative limit which one can identify with the algebra of
rational functions on the 2-plane minus the axes.  We shall see that the
commutation relations between the elements of the algebra and the
1-forms determine a set of torsion-free metric connections which remain
non-trivial in the commutative limit. It is to be expected that the
converse is true, that every torsion-free metric connection on the
2-plane determines a set of commutation relations. This would mean that
in particular the covariant calculus of Wess \& Zumino (1990) is
determined by a geometry on the 2-plane.

The differential calculi we introduce here are based on derivations and
are similar to those introduced by Dubois-Violette (1988) and
Dubois-Violette {\it et al.} (1989) to construct differential calculi
over matrix algebras. We refer to Madore (1995a,b) for a detailed
description of how to use this method to construct a sequence of
differential calculi over a given matrix algebra.  In this previous work
the set of derivations were chosen to form the Lie algebra of the
special linear group $SL_m$. With this restriction the number of
differential calculi which can be put on a matrix algebra of dimension
$n$ is equal to the number of integers $m$ such that the $SL_m$ has an
irreducible representation on a space of dimension $n$. There are always
of course at least two, $m = 2$ and $m = n$ but for large $n$ there are
in general many more.

In Section~2 we present a general method for constructing differential
calculi, based on a set of derivations which do not necessarily close to
form a Lie algebra. In Sections~3 and~4 we present some examples using
as algebra the generalized quantum plane. In Section~5 we investigate
linear connections and show how they depend on the choice of
differential calculus as well as, of course, the algebra. By ${\cal A}$
we designate an arbitrary associative algebra with unit element and with
center ${\cal Z}({\cal A})$

\section{General Formalism}

Of the many differential calculi which can be constructed over an
algebra ${\cal A}$ the largest is the differential envelope or universal
differential calculus $(\Omega^*_u({\cal A}), d_u)$. Every other
differential calculus can be considered as a quotient of it.  For the
definitions we refer, for example, to the book by Connes (1994). Let
$(\Omega^*({\cal A}), d)$ be another differential calculus over 
${\cal A}$. Then there exists a unique $d_u$-homomorphism $\phi$
$$
\def\normalbaselines{\baselineskip=18pt}
\matrix{
{\cal A} &\buildrel d_u \over \longrightarrow &\Omega_u^1({\cal A})
         &\buildrel d_u \over \longrightarrow &\Omega_u^2({\cal A})   
         &\buildrel d_u \over \longrightarrow &\cdots                \cr
\parallel&&\phi_1 \downarrow \phantom{\phi_1}
         &&\phi_2 \downarrow \phantom{\phi_2}                        \cr
{\cal A} &\buildrel d \over \longrightarrow &\Omega^1({\cal A})
         &\buildrel d \over \longrightarrow &\Omega^2({\cal A})   
         &\buildrel d \over \longrightarrow &\cdots                  \cr
}
\def\normalbaselines{\baselineskip=12pt}                           \eqno(2.1)
$$
of $\Omega^*_u({\cal A})$ onto $\Omega^*({\cal A})$. It is given by
$$
\phi (d_u a) =  d a.                                               \eqno(2.2)
$$
The restriction $\phi_p$ of $\phi$ to each $\Omega^p_u$ is defined by
$$
\phi_p(a_0 d_u a_1 \cdots d_u a_p) = a_0 da_1 \cdots da_p.
$$

Consider a given algebra ${\cal A}$ and suppose that we know how to
construct an ${\cal A}$-module $\Omega^1({\cal A})$ and an application
$$
{\cal A} \buildrel d \over \longrightarrow \Omega^1({\cal A}).     \eqno(2.3)
$$
Then using (2.1) there is a method of constructing $\Omega^p({\cal A})$
for $p \geq 2$ as well as the extension of the differential.  Since we
know $\Omega_u^1({\cal A})$ and $\Omega^1({\cal A})$ we can suppose that
$\phi_1$ is given. We must construct $\Omega^2({\cal A})$. The simplest
consistent choice would be to set
$$
\Omega^2({\cal A}) = \Omega_u^2({\cal A})/d_u {\rm Ker} \phi_1.    \eqno(2.4)
$$
This is the largest differential calculus consistent with the
constraints on $\Omega^1({\cal A})$.  The map $\phi_2$ is defined to be
the projection of $\Omega_u^2({\cal A})$ onto $\Omega^2({\cal A})$ so
defined and $d$ is defined by $d(fdg) = df dg$.  This procedure can be
continued by iteration to arbitrary order in $p$. See, for example, Baehr
{\it et al.} (1995).

To initiate the above construction we shall define the 1-forms using a
set of derivations.  We shall suppose that they are interior and exclude
therefore the case where ${\cal A}$ is commutative.  For each integer
$n$ let $\lambda_a$ be a set of $n$ linearly independent elements of
${\cal A}$ and introduce the derivations $e_a = {\rm ad}\, \lambda_a$.
In general the $e_a$ do not form a Lie algebra but they do however
satisfy commutation relations as a consequence of the commutation
relations of ${\cal A}$. Define the map (2.3) by
$$
df (e_a) = e_a \, f.                                               \eqno(2.5)
$$
We shall suppose that there exists a set of $n$ elements $\theta^a$ of
$\Omega^1({\cal A})$ such that
$$
\theta^a (e_b) = \delta^a_b.                                       \eqno(2.6)
$$
In the examples which we consider we shall show that the $\theta^a$
exist by explicit construction.  We shall refer to the set of $\theta^a$
as a frame or Stehbein.  It commutes with the elements $f \in {\cal A}$,
$$
f \theta^a = \theta^a f.                                          \eqno(2.7)
$$

The construction of the $\theta^a$ from the derivations did not really
use the fact that they were inner.  For example if the $e_a$ are $n$
linearly independent vector fields on a manifold $V$ of dimension $n$,
that is, $n$ linearly independent outer derivations of the algebra
${\cal A} = {\cal C}(V)$ of smooth functions on $V$ then 
$\Omega^*({\cal A})$ is the algebra of de~Rham forms.

The ${\cal A}$-bimodule $\Omega^1({\cal A})$ is generated by all
elements of the form $fdg$ or of the form $(df)g$. Because of the
Leibniz rule these conditions are equivalent.  Define 
$\theta = - \lambda_a \theta^a$. Then one sees that
$$
df = - [\theta, f]                                                 \eqno(2.8)
$$
and it follows that as a bimodule $\Omega^1({\cal A})$ is generated by
one element.  It follows also that the 2-form $d\theta + \theta^2$ can
be written in the form
$$
d\theta + \theta^2 = - {1\over 2} K_{ab} \theta^a \theta^b         \eqno(2.9)
$$
with coefficients $K_{ab}$ which lie in ${\cal Z}({\cal A})$. By 
definition
$$
f dg (e_a) = f e_a \, g, \qquad (dg) f (e_a) =  (e_a \, g) f.
$$
Using the frame we can write these as
$$
f dg = (f e_a g) \theta^a, \qquad (dg) f =  (e_a g) f \theta^a.    \eqno(2.10)
$$
The commutation relations of the algebra constrain the relations
between $f dg$ and $(dg) f$ for all $f$ and $g$. 

As a left or right module, $\Omega^1({\cal A})$ is free of rank $n$.
Because of the commutation relations of the algebra or, equivalently,
because of the kernel of $\phi_1$ in the quotient (2.4) the $\theta^a$
satisfy commutation relations. We shall suppose that they can be written
in the form
$$
\theta^a \theta^b + C^{ab}{}_{cd} \theta^c \theta^d = 0.           \eqno(2.11)
$$
If $C^{ab}{}_{cd} = \delta^a_c\delta^b_d$ then one sees that
$\Omega^2({\cal A}) = 0$.  It follows from (2.11) that for an arbitrary
element $f$ of the algebra
$$
[f, C^{ab}{}_{cd}] \theta^c \theta^d = 0.
$$
We shall suppose that
$$
C^{ab}{}_{ef} C^{ef}{}_{cd} = \delta^a_c \delta^b_d
$$
and that the relations (2.11) are complete in the sense that if
$A_{ab} \theta^a \theta^b = 0$ we can conclude that
$$
A_{ab} - C^{cd}{}_{ab} A_{cd} = 0.                                 \eqno(2.12)
$$
This will be the case for all the differential calculi which we shall
consider on the generalized quantum plane in the next sections.  We can
conclude then that the $C^{ab}{}_{cd}$ are elements of ${\cal Z}({\cal A})$. 
In ordinary geometry one can choose
$$
C^{ab}{}_{cd} = \delta^b_c \delta^a_d 
$$
and the relation (2.11) expresses the fact that the 1-forms anticommute.
Let $\bigwedge_C^*$ be the twisted exterior algebra determined by the
the relations (2.11). Then the differential algebra $\Omega^*({\cal A})$
can be factorized as the tensor product
$$
\Omega^*({\cal A}) = {\cal A} \otimes_\bbbc {\textstyle \bigwedge}_C^*.
$$

Because the 2-forms are generated by products of the $\theta^a$ one has
$$
d\theta^a = - {1\over 2} C^a{}_{bc} \theta^b \theta^c.             \eqno(2.13) 
$$
Without loss of generality we can suppose that the structure elements
$C^a{}_{bc}$ satisfy the identities
$$
C^a{}_{bc} + C^a{}_{de} C^{de}{}_{bc} = 0.                         \eqno(2.14)
$$
It is to be noticed that they do not in general lie ${\cal Z}({\cal A})$
.
In fact from the identity $d(f\theta^a) = d(\theta^a f)$ one
sees that
$$
\Big({1\over 2}[C^a{}_{bc}, f] +
e^{\phantom a}_{(b} f \delta^a_{c)}\Big) \theta^b \theta^c = 0.    \eqno(2.15)
$$
Using the definition of the derivations one can write this in the form
$$
\Big({1\over 2} C^a{}_{bc} + \lambda^{\phantom a}_{(b} \delta^a_{c)} 
- {1\over 2} D^a{}_{bc} \Big)  \theta^b \theta^c = 0               \eqno(2.16)
$$
with $D^a{}_{bc} \in {\cal Z}({\cal A})$. We can suppose that the 
$D^a{}_{bc}$ satisfy (2.14):
$$
D^a{}_{bc} + D^a{}_{de} C^{de}{}_{bc} = 0.                         \eqno(2.17)
$$

Using this and the relations (2.11) and (2.14) as well as the
completeness assumption (2.12) we can conclude from (2.16) that
$$
C^a{}_{bc} - D^a{}_{bc} + \lambda^{\phantom a}_{(b} \delta^a_{c)} 
- \lambda^{\phantom a}_{(d} \delta^a_{e)} C^{de}{}_{bc} = 0.       \eqno(2.18)
$$

The equation (2.16) can also be written in the form
$$
d\theta^a = - [\theta, \theta^a] 
            - {1\over 2} D^a{}_{bc} \theta^b \theta^c              \eqno(2.19)
$$
with a graded commutator. If $D^a{}_{bc} = 0$ the form (2.8) for the
exterior derivative is valid for all elements of $\Omega^*({\cal A})$
and the element $\theta$ plays the role of the phase of a generalized
Dirac operator in the sense of Connes (1986).

From (2.19) we find that
$$
d\theta = - 2 \theta^2 
+ {1 \over 2} \lambda_a D^a{}_{bc} \theta^b \theta^c.
$$
Comparing this with (2.9) we conclude that 
$$
\theta^2 =  {1 \over 2}(\lambda_a D^a{}_{bc} +
                                      K_{bc}) \theta^b \theta^c.   \eqno(2.20)
$$
If we suppose that $K_{bc}$ satisfies (2.14), 
$$
K_{ab} + K_{cd} C^{cd}{}_{ab} = 0,
$$
then we can conclude that
$$
(\lambda_b \lambda_c - C^{de}{}_{bc} \lambda_d \lambda_e
- \lambda_a D^a{}_{bc} - K_{bc}) \theta^b \theta^c = 0.            \eqno(2.21)
$$
Using again the completeness assumption (2.12) we find
$$
\lambda_b \lambda_c - C^{de}{}_{bc} \lambda_d \lambda_e = 
\lambda_a D^a{}_{bc} + K_{bc}.                                     \eqno(2.22)
$$
If we introduce the twisted bracket
$$
[\lambda_a, \lambda_b]_C = 
 \lambda_a  \lambda_b - C^{cd}{}_{ab} \lambda_c \lambda_d. 
$$
we can rewrite (2.22) in the form
$$
[\lambda_b, \lambda_c]_C = \lambda_a D^a{}_{bc} + K_{bc}.          \eqno(2.23)
$$
If we write out in detail the equation $d^2 f = 0$, using (2.12) -
(2.14) we find that
$$
[e_b, e_c]_C f = e_a f C^a{}_{bc}.                                 \eqno(2.24)
$$
This is the dual relation to the `Maurer-Cartan' equation (2.13).

The constraint (2.23) follows from the commutation relations (2.11) we
have supposed for the frame as well as from the conditions we have
imposed on the coefficients $C^{ab}{}_{cd}$. In the matrix case the
general formalism simplifies considerably. The $\theta^a$ are given in
terms of $d\lambda^a$ by
$$
\theta^a = \lambda_b \lambda^a d\lambda^b.
$$
The elements of the frame anticommute and one can choose 
$C^{ab}{}_{cd} = \delta^b_c \delta^a_d$.  In Equation~(2.19) the first
term on the right-hand side vanishes and $D^a{}_{bc} = C^a{}_{bc}$. On
the right-hand side of Equations~(2.9) and (2.23) we have $K_{ab} = 0$.

\section{Calculi based on 2 derivations}

Using the construction of the previous section one can construct an
infinite sequence of differential calculi over the generalized quantum
plane ${\cal A}$, the algebra generated by four elements 
$(x, y, x^{-1}, y^{-1})$ subject to the relation $xy = q yx$ as well as
the usual relations between an element and its inverse. Here $q$ is an
arbitrary complex number. The subalgebra generated by $(x, y)$ alone
with the covariant differential calculus of Wess \& Zumino (1990) is
called the quantum plane.  In this section we shall consider the case 
$n = 2$.  Define, for $q \neq 1$
$$
\lambda_1 = {q \over q - 1} y,                \qquad
\lambda_2 = {q \over q - 1} x.                                    \eqno(3.1)
$$
The normalization has been chosen here so that the structure elements
$C^a{}_{12}$ contain no factors $q$. The $\lambda_a$ are singular in the
limit $q \rightarrow 1$ for the same reason as the limit from quantum 
mechanics:
$$
{1\over \hbar} {\rm ad}\, p \rightarrow 
{1\over i} {\partial\over \partial q}.
$$
We find that
$$
e_1 x = - x y,   \qquad
e_1 y = 0,       \qquad
e_2 x = 0,       \qquad
e_2 y =   x y.                                                    \eqno(3.2)
$$
These rather unusual derivations are extended to arbitrary
polynomials in the generators by the Leibniz rule.  From (3.2) we
conclude that the commutation relations which follow from (2.10) are
$$
\matrix{
x dx = q dx x,        \hfill   &\qquad&  
y dx = q^{-1} dx y,   \hfill\cr&&\cr
x dy = q dy x,        \hfill   &\qquad&
y dy = q^{-1} dy y.   \hfill
}                                                                  \eqno(3.3)
$$
From these relations if $q \neq -1$ we deduce
$$
(dx)^2 = 0, \qquad (dy)^2 = 0,\qquad  dx dy + q dy dx = 0.         \eqno(3.4)
$$
Using the relations (3.2) we find
$$
d x =  -  x y \theta^1,                             \qquad
d y =     x y \theta^2                                             \eqno(3.5)
$$
and solving for the $\theta^a$ we obtain
$$
\theta^1 = - q^{-1} x^{-1} y^{-1} dx,               \qquad
\theta^2 =   q^{-1} x^{-1} y^{-1} dy.                              \eqno(3.6) 
$$
The $\theta^a$ satisfy the commutation relations
$$
(\theta^1)^2=0,                                     \qquad 
(\theta^2)^2=0,                                     \qquad
\theta^1\theta^2 + q \theta^2\theta^1 = 0.                         \eqno(3.7)
$$
This is of the form (2.11). If we reorder the indices 
$(11, 12, 21, 22) = (1,2,3,4)$ then the matrix $C^{ab}{}_{cd}$ can be
written in the form of a $4 \times 4$ matrix
$$
C = \pmatrix{
1 &0      &0 &0 \cr
0 &0      &q &0 \cr
0 &q^{-1} &0 &0 \cr
0 &0      &0 &1
}.                                                                 \eqno(3.8)
$$
That is, $C^{12}{}_{21} = q$ and $C^{21}{}_{12} = q^{-1}$.
The  structure elements $C^a{}_{bc}$ are given by
$$
C^1{}_{12} = - x,                                    \qquad 
C^2{}_{12} = - y                                                   \eqno(3.9) 
$$
and the condition (2.14). Equation (2.18) is satisfied.  For
$\theta$ we find the expression
$$
\theta = {1 \over q - 1} (q x^{-1} dx - y^{-1} dy).               \eqno(3.10)
$$
It is a closed form.

As a second example with $n=2$ we define, for $q^4 \neq 1$
$$
\lambda_1 = {1 \over q^4 - 1} x^{-2} y^2,             \qquad
\lambda_2 = {1 \over q^4 - 1} x^{-2}.                             \eqno(3.11)
$$
The normalization has been chosen here again so that the structure
elements $C^a{}_{12}$ contain no factors $q$. We find then that for 
$q^2 \neq -1$
$$
\matrix{
\displaystyle{e_1 x = - {1 \over q^2 (q^2 + 1)}  x^{-1} y^2}, \hfill 
&\displaystyle{e_1 y = - {1 \over q^2 + 1} x^{-2} y^3},       \hfill\cr   
\displaystyle{e_2 x = 0},                                     \hfill
&\displaystyle{e_2 y = - {1 \over q^2 + 1} x^{-2} y}.         \hfill 
}                                                                 \eqno(3.12)
$$
From these we conclude that the commutation relations which follow from
(2.10) are
$$
\matrix{
xdx = q^2 dx x,               \hfill   &\qquad&  
x dy = q dy x + (q^2-1)dx y,  \hfill\cr&&\cr
ydx = q dx y,                 \hfill   &\qquad&
y dy = q^2 dy y.              \hfill
}                                                                 \eqno(3.13)
$$
We have then in this case the covariant differential calculus of Wess \&
Zumino (1990).  It has been encoded in the functional form of the
$\lambda_a$.  If $q^2 \neq -1$ from (3.13) we deduce
$$
(dx)^2 = 0, \qquad (dy)^2 = 0,\qquad  dy dx + q dx dy = 0.        \eqno(3.14)
$$
Using the relation (2.6) we find
$$
d x = - {1 \over q^2 (q^2 + 1)} x^{-1} y^2 \theta^1, \qquad
d y = - {1 \over q^2 + 1} x^{-2} y (y^2 \theta^1 + \theta^2)       \eqno(3.15)
$$
and solving for the $\theta^a$ we obtain
$$
\theta^1 = - q^4 (q^2 + 1) x y^{-2} dx,              \qquad
\theta^2 = - q^2 (q^2 + 1) x (x y^{-1} dy - dx).                   \eqno(3.16) 
$$
The $\theta^a$ satisfy the commutation relations
$$
(\theta^1)^2=0,                                       \qquad 
(\theta^2)^2=0,                                       \qquad 
q^4 \theta^1\theta^2 + \theta^2\theta^1 = 0.                       \eqno(3.17)
$$
This is of the form (2.11) if the matrix $C^{ab}{}_{cd}$ is
given by the $4 \times 4$ matrix
$$
C = \pmatrix{
1 &0   &0      &0 \cr
0 &0   &q^{-4} &0 \cr
0 &q^4 &0      &0 \cr
0 &0   &0      &1
}.                                                                 \eqno(3.18)
$$
That is, $C^{12}{}_{21} = q^{-4}$ and $C^{21}{}_{12} = q^4$.
The structure elements $C^a{}_{bc}$ are given by
$$
C^1{}_{12} = - x^{-2},                                  \qquad 
C^2{}_{12} = - x^{-2} y^2                                          \eqno(3.19) 
$$
and the condition (2.14). Equation (2.18) is again satisfied.

For $\theta$ we find the expression
$$
\theta = {q^2 \over q^2 - 1} y^{-1} dy.                            \eqno(3.20)
$$
It is again a closed form.
 
From these two examples we see that each choice of two elements
$\lambda_1$ and $\lambda_2$ gives rise to a differential calculus on the
generalized quantum plane and we shall see in the Section~5 that each 
choice gives rise to a linear connection.

\section{Calculi based on 3 derivations} %4

In this section we shall consider the case $n = 3$. There is an
essential difference with the previous case in that relations of the
form (3.3) or (3.13) which allow one to pass from one side of the
differential to the other no longer hold. The difference is given in
fact in terms of the extra elements of the frame. What we do is
extend the definition of $dx$ and $dy$ to another derivation and the
extension satisfies quite naturally less relations. The left (or right)
module $\Omega^1({\cal A})$ is now of rank 3 instead of 2. As an example
we extend the $\lambda_a$ defined in (3.1) by the addition of a
$\lambda_3$:
$$
\lambda_1 = {q \over q - 1} y,                \qquad
\lambda_2 = {q \over q - 1} x,                \qquad 
\lambda_3 = {q \over q - 1} \alpha x y.                            \eqno(4.1)
$$
The $\alpha$ is an arbitrary complex number.  We have then 
$[\lambda_1, \lambda_2] = - \alpha^{-1} \lambda_3$ but of course the set
of $\lambda_a$ do not form a Lie algebra. To the relations (3.2) we must
add two additional ones,
$$
e_3 x = - \alpha x^2 y, \qquad e_3 y = \alpha x y^2,               \eqno(4.2)
$$
and so we find
$$
d x =  -  x y \theta^1 - \alpha x^2 y \theta^3,   \qquad
d y =     x y \theta^2 + \alpha x y^2 \theta^3.                    \eqno(4.3)
$$
instead of (3.5). Define
$$
\tau = x dy - q dy x.                                               \eqno(4.4)
$$
Then one of the commutation relations (3.3) becomes an expression for 
$\theta^3$ in terms of $\tau$:
$$
\tau = \alpha q^{-1} (q - 1) x^2 y^2 \theta^3.                      \eqno(4.5)
$$
We can solve then (4.3) for the $\theta^a$ and we obtain
$$
\eqalign{
&\theta^1 = - q^{-1} x^{-1} y^{-1} dx 
           - {1\over q^2(q-1)} x^{-1} y^{-2} \tau,          \cr
&\theta^2 =   q^{-1} x^{-1} y^{-1} dy
           - {1\over q(q-1)} x^{-2} y^{-1} \tau,            \cr
&\theta^3 = {1\over \alpha q^3(q-1)} x^{-2} y^{-2} \tau 
}                                                                  \eqno(4.6) 
$$          
instead of (3.6).  This frame is singular in the limit $q \rightarrow 1$ 
as it must be. The differential calculus, expressed in terms of $dx$ and
$dy$, has however a well-defined limit which lies somewhere between the
de~Rham calculus and the universal one. For a discussion of this point
we refer to Dimakis \& M\"uller-Hoissen (1992), Dimakis \& Tzanakis
(1995) and to Baehr {\it et al.} (1995).

If $q \neq -1$ we can deduce the first two of the relations (3.4) and we
can conclude that
$$
\eqalign{
&q (\theta^1)^2 + \alpha x (\theta^1 \theta^3 + q \theta^3 \theta^1)
               + \alpha^2 x^2 (\theta^3)^2 = 0,                    \cr
&q (\theta^2)^2 + \alpha y (\theta^3 \theta^2 + q \theta^2 \theta^3)
               + \alpha^2 y^2 (\theta^3)^2 = 0.
}                                                                   \eqno(4.7)
$$
Multiply the first equation by $y$ and the second by $x$ and commute
through. One deduces then that each of the coefficients must vanish:
$$
(\theta^1)^2 = 0, \qquad (\theta^2)^2 = 0,  \qquad (\theta^3)^2 = 0,
$$
and
$$
\theta^1 \theta^3 + q \theta^3 \theta^1 = 0,      \qquad
\theta^3 \theta^2 + q \theta^2 \theta^3 = 0.                       \eqno(4.8)
$$
There is missing a relation between $\theta^1\theta^2$ and
$\theta^2\theta^1$. We must therefore rather artificially complete the
coefficients in (2.11) by setting $C^{12}{}_{12} = -1$ and
$C^{12}{}_{21} = 0$. From (2.23) we find then that $K_{ab}=0$ and the
$D^a{}_{bc}$ vanish except for the values
$$
D^3{}_{12} = {2\over\alpha(q-1)}, \qquad
D^3{}_{21} = q D^3{}_{12}.                                         \eqno(4.9)
$$

The form $\theta$ is given by
$$
\theta = - {q \over q-1} (y \theta^1 + x \theta^2 + \alpha x y \theta^3).
$$
It follows then immediatedly from (2.19) that
$$
\eqalign{
&d\theta^1  =  
{q\over q-1} x (\theta^1\theta^2+\theta^2\theta^1) +
\alpha xy \theta^1\theta^3,                              \cr
&d\theta^2  =  
{q\over q-1} y (\theta^1\theta^2+\theta^2\theta^1) +
\alpha xy \theta^3\theta^2,                              \cr
&d\theta^3  = y \theta^1\theta^3 + x \theta^3\theta^2 
- {1\over\alpha(q-1)}(\theta^1\theta^2+q\theta^2\theta^1).
}                                                                 \eqno(4.10)
$$
from which one can read off the expressions for the structure elements
which extend (3.9).  The third of the relations (3.4) becomes
$$
d\tau = - x^2 y^2\Big((\theta^1 \theta^2 + q \theta^2 \theta^1)
        + \alpha x (\theta^2 \theta^3 + \theta^3 \theta^2)
        + \alpha y (\theta^3 \theta^1 + \theta^1 \theta^3)\Big).
$$
Using (4.10) one finds
$$
d\tau = - x^2 y^2\Big((\theta^1 \theta^2 + q \theta^2 \theta^1)
        - \alpha q^{-1} (q - 1) d \theta^3 \Big).                 \eqno(4.11)
$$

If one adds to (4.8) the supplementary relation
$$
\theta^1 \theta^2 + q \theta^2 \theta^1 = 0                       \eqno(4.12)
$$
then $\Omega^2({\cal A})$ becomes a quotient of the right-hand side of
(2.4). We can set $C^{12}{}_{21} = q$ and $C^{12}{}_{12} = 0$ as in
(3.8).  Now we have $K_{ab} = 0$ and $D^a{}_{bc} = 0$ and
Equations~(4.10) reduce to
$$
\eqalign{
&d\theta^1 = x \theta^1\theta^2 + \alpha xy \theta^1\theta^3,    \cr
&d\theta^2 = y \theta^1\theta^2 + \alpha xy \theta^3\theta^2,    \cr
&d\theta^3 = y \theta^1\theta^3 + x \theta^3\theta^2.
}
$$
Equation~(4.11) simplifies to
$$
d\tau = \alpha x^2 y^2 q^{-1} (q - 1) d \theta^3.                  \eqno(4.13)
$$

A similar extension of the second example of the previous section can be
given, again by introducing a third derivation. As before this yields an
extension of the module of forms as a left (or right) module.

\section{Linear connections} %5

For each of the differential calculi defined in the previous section one
can define a set of linear connections.  The definition of a connection
as a covariant derivative was given an algebraic form in the Tata
lectures by Koszul (1960) and generalized to noncommutative geometry by
Karoubi (1981) and Connes (1986, 1994). We shall often use here the
expressions `connection' and `covariant derivative' synonymously.  In
fact we shall distinguish three different types of connections. A `left
${\cal A}$-connection' is a connection on a left ${\cal A}$-module; it
satisfies a left Leibniz rule.  A `bimodule ${\cal A}$-connection' is a
connection on a general bimodule ${\cal M}$ which satisfies a left and
right Leibniz rule. In the particular case where ${\cal M}$ is the
module of 1-forms we shall speak of a `linear connection'. A bimodule
over an algebra ${\cal A}$ is also a left module over the tensor product
${\cal A}^e = {\cal A} \otimes_\bbbc {\cal A}^{\rm op}$ of the algebra
with its `opposite'. Such a bimodule can have a bimodule 
${\cal A}$-connection as well as a left ${\cal A}^e$-connection. (Cuntz
\& Quillen 1995, Bresser {\it et al.} 1995).  These two definitions are
compared in Dubois-Violette {\it et al.} (1995b).

Let ${\cal A}$ be an arbitrary algebra and $(\Omega^*({\cal A}) ,d)$ a
differential calculus over ${\cal A}$. One defines a left 
${\cal A}$-connection on a left ${\cal A}$-module ${\cal H}$ as a
covariant derivative
$$
{\cal H} \buildrel D \over \rightarrow 
\Omega^1({\cal A}) \otimes_{\cal A} {\cal H}                      \eqno(5.1) 
$$ 
which satisfies the left Leibniz rule
$$
D (f \psi) =  df \otimes \psi + f D\psi                           \eqno(5.2)
$$
for arbitrary $f \in {\cal A}$. This map has a natural extension
$$
\Omega^*({\cal A}) \otimes_{\cal A} {\cal H} 
\buildrel \nabla \over \longrightarrow
\Omega^*({\cal A}) \otimes_{\cal A} {\cal H}                      \eqno(5.3)
$$
given, for $\psi \in {\cal H}$ and $\alpha \in \Omega^n({\cal A})$,
by $\nabla \psi = D \psi$ and 
$$
\nabla (\alpha \psi) = d\alpha \otimes \psi + 
                (-1)^n \alpha \nabla \psi.
$$
The operator $\nabla^2$ is necessarily left-linear. However when 
${\cal H}$ is a bimodule it is not in general right-linear.

A covariant derivative on the module $\Omega^1({\cal A})$ must satisfy
(5.2). But $\Omega^1({\cal A})$ has also a natural structure as a right
${\cal A}$-module and one must be able to write a corresponding right
Leibniz rule in order to construct a bilinear curvature. Quite generally
let ${\cal M}$ be an arbitrary bimodule. Consider a covariant derivative
$$
{\cal M} \buildrel D \over \rightarrow 
\Omega^1({\cal A}) \otimes_{\cal A} {\cal M}                       \eqno(5.4) 
$$ 
which satisfies both a left and a right Leibniz rule.  In order to
define a right Leibniz rule which is consistent with the left one, it
was proposed by Mourad (1995), by Dubois-Violette \& Michor (1995) and
by Dubois-Violette \& Masson (1995) to introduce a generalized
permutation
$$
{\cal M} \otimes_{\cal A} \Omega^1({\cal A})
\buildrel \sigma \over \longrightarrow
\Omega^1({\cal A}) \otimes_{\cal A} {\cal M}.
$$
The right Leibniz rule is given then as
$$
D(\xi f) = \sigma (\xi \otimes df) + (D\xi) f                     \eqno(5.5)
$$
for arbitrary $f \in {\cal A}$ and $\xi \in {\cal M}$. The purpose of
the map $\sigma$ is to bring the differential to the left while
respecting the order of the factors. It is necessarily bilinear
(Dubois-Violette {\it et al.} 1995a). Let $\pi$ be the product in the
algebra of forms.  It was argued by Mourad (1995) and by Dubois-Violette
{\it et al.} (1995a) that a necessary as well as sufficient condition
for torsion to be right-linear is that $\sigma$ satisfy the consistency
condition
$$
\pi \circ (\sigma + 1) = 0.                                        \eqno(5.6) 
$$
We define a bimodule ${\cal A}$-connection to be the couple
$(D, \sigma)$. We shall make no mention of curvature. There is at the
moment no general concensus of the correct definition of the curvature 
of a bimodule connection. The problem is that the operator $\nabla^2$ 
is not in general right-linear.  We refer to Dubois-Violette {\it et
al.} (1995b) for a recent discussion.

This general formalism can be applied in particular to the differential
calculi which we have constructed in Section~2.  Since $\Omega^1({\cal A})$ 
is a free module the map $\sigma$ can be defined by the action on the
basis elements:
$$
\sigma (\theta^a \otimes \theta^b) = 
S^{ab}{}_{cd} \theta^c \otimes \theta^d.                            \eqno(5.7)
$$
By the sequence of identities
$$
f S^{ab}{}_{cd} \theta^c \otimes \theta^d = 
\sigma (f \theta^a \otimes \theta^b) = 
\sigma (\theta^a \otimes \theta^b f) =
S^{ab}{}_{cd} f \theta^c \otimes \theta^d
$$
we conclude that the coefficients $S^{ab}{}_{cd}$ must lie in 
${\cal Z}({\cal A})$.  From (2.12) we see that the condition (5.6) can
be written
$$
(\delta^a_e\delta^b_f + S^{ab}{}_{ef})
(\delta^e_c\delta^f_d - C^{ef}{}_{cd}) = 0.                         \eqno(5.8)
$$
A natural, but certainly not the unique, choice of $\sigma$ is
given by $S^{ab}{}_{cd} = C^{ab}{}_{cd}$.

Since $\Omega^1({\cal A})$ is a free module a covariant
derivative can be defined by its action on the basis elements:
$$
D\theta^a = - \omega^a{}_{bc} \theta^b \otimes \theta^c.            \eqno(5.9)
$$
The coefficients here are elements of the algebra.  Because of the
identity $D(f \theta^a) = D(\theta^a f)$ there are very stringent
compatibility conditions, which using (5.7) can be written out as
$$
[\omega^a{}_{bc}, f] = 
       e_d f (S^{ad}{}_{bc} - \delta^d_b \delta^a_c).
$$
The general solution to this equation is of the form
$\omega^a{}_{bc} = \omega_{(0)}{}^a{}_{bc} + \chi^a{}_{bc}$ where 
$$
\omega_{(0)}{}^a{}_{bc} = 
    \lambda_d (S^{ad}{}_{bc} - \delta^d_b \delta^a_c)               \eqno(5.10)
$$
and $\chi^a{}_{bc} \in {\cal Z}({\cal A})$.  One can also express
$D_{(0)}$ in the form (Dubois-Violette {\it et al.} 1995a, Madore {\it
et al.} 1995)
$$
D_{(0)} \theta^a = - \theta \otimes \theta^a + 
            \sigma (\theta^a \otimes \theta).
$$
 
The torsion 2-form is defined as usual as
$$
\Theta^a = d \theta^a - \pi \circ D \theta^a
$$
Comparing (5.10) with (2.19), we see that the torsion $\Theta_{(0)}^a$
of $D_{(0)}$ is given by
$$
\Theta_{(0)}^a = - {1\over 2} D^a{}_{bc} \theta^b\theta^c.         \eqno(5.11)
$$
In general a covariant derivative is torsion-free provided the condition
$$
\omega^a{}_{bc} - \omega^a{}_{de} C^{de}{}_{bc} = C^a{}_{bc}
$$
is satisfied. The covariant derivative (5.9) is torsion free if and only
if
$$
\chi^a{}_{bc} = {1\over 2} D^a{}_{bc}.
$$

On the ordinary quantum plane one can show that there is there is a
unique 1-parameter family of linear connections (Dubois-Violette {\it et
al.} 1995a) and that this connection is torsion free.  We find here a
different result; there is an ambiguity which depends on elements of
${\cal Z}({\cal A})$. An interesting limit case is given by
$$
S^{ab}{}_{cd} = C^{ab}{}_{cd} = \delta^b_c \delta^a_d.             \eqno(5.12)
$$
In this case from (2.18) one sees that $D^a{}_{bc} = C^a{}_{bc} \neq 0$.
From (2.22) one sees that $K_{ab} = 0$ and the $\lambda_a$ form a Lie
algebra.  The matrix case is a typical example.  From (5.10) if follows
that $D_{(0)} \theta^a = 0$ and so $D_{(0)}$ has torsion but no
curvature. The connection corresponds to the canonical flat connection
on a parallelizable manifold.

One can define a metric by the condition
$$
g(\theta^a \otimes \theta^b) = g^{ab}                              \eqno(5.13)
$$
where the coefficients $g^{ab}$ are elements of the algebra. To be well
defined on all elements of the tensor product 
$\Omega^1({\cal A}) \otimes_{\cal A} \Omega^1({\cal A})$ the metric must
be bilinear and by the sequence of identities
$$
f g^{ab} = g(f \theta^a \otimes \theta^b) 
= g(\theta^a \otimes \theta^b f) = g^{ab} f                      \eqno(5.14)
$$
we conclude that the coefficients must lie in ${\cal Z}({\cal A})$.  The
covariant derivative (5.9) is compatible with the metric
(Dubois-Violette {\it et al.} 1995a) if and only if
$$
\omega^a{}_{bc} + \omega_{ce}{}^f S^{ae}{}_{bf} = 0.              \eqno(5.15)
$$
The condition that (5.10) be metric compatible can be written as
$$
S^{ae}{}_{dh} g^{hf} S^{cb}{}_{ef} = g^{ac} \delta^b_d.           \eqno(5.16)
$$

Consider now the first differential calculus of Section~3, defined by
(3.1).  On the right-hand side of (2.23) we have $K_{ab} = 0$ and 
$D^a{}_{bc} = 0$. The torsion of $D_{(0)}$ vanishes.  The coefficients
$g^{ab}$ are complex numbers.  With the convention of (3.8) they can be
written as $(g^1, g^2, g^3, g^4)$. Using the $GL(2, \bbbc)$-invariance
one can impose that
$$
g^4 = g^1, \qquad g^3 = - g^2.
$$ 
If we suppose that $g^2 = 0$ there is no restriction in supposing that
$g^1 = 1$; the $g^{ab}$ are the components of the euclidean metric in
two dimensions.  With the convention of (3.8) the condition (5.16) can
be written in the matrix form
$$ 
\pmatrix{
S^1{}_1 & S^1{}_2 & S^1{}_3 & S^1{}_4  \cr
S^2{}_1 & S^2{}_2 & S^2{}_3 & S^2{}_4  \cr
S^3{}_1 & S^3{}_2 & S^3{}_3 & S^3{}_4  \cr
S^4{}_1 & S^4{}_2 & S^4{}_3 & S^4{}_4
}
\pmatrix{
S^1{}_1 & S^1{}_3 & S^3{}_1 & S^3{}_3  \cr
S^1{}_2 & S^1{}_4 & S^3{}_2 & S^3{}_4  \cr
S^2{}_1 & S^2{}_3 & S^4{}_1 & S^4{}_3  \cr
S^2{}_2 & S^2{}_4 & S^4{}_2 & S^4{}_4
} = 1.                                                             \eqno(5.17)
$$
From the approximation linear in $q - 1$ one sees that the solution must
be of the form
$$ 
S = \pmatrix{
S^1{}_1 & 0       & 0       & S^1{}_4  \cr
0       & S^2{}_2 & S^2{}_3 & 0        \cr
0       & S^3{}_2 & S^3{}_3 & 0        \cr
S^4{}_1 & 0       & 0       & S^4{}_4
}.                                                                 \eqno(5.18)
$$
The consistancy conditions (5.8) become
$$
1 + S^2{}_2 = q^{-1} S^2{}_3, \qquad 
1 + S^3{}_3 = q S^3{}_2.                                           \eqno(5.19)
$$
In general $S^{ab}{}_{cd} = C^{ab}{}_{cd}$ does not yield a
metric-compatible covariant derivative.  There is a solution however to
(5.17), (5.19) given by
$$ 
S = {1 \over q^2 + 1}\pmatrix{
2 q     & 0       & 0       & 1 - q^2   \cr
0       & 1 - q^2 & 2 q     & 0         \cr
0       & 2 q     & q^2 - 1 & 0         \cr
q^2 - 1 & 0       & 0       & 2 q
}                                                                 \eqno(5.20)
$$
That is, for example, 
$$
S^{12}{}_{21} = S^{21}{}_{12} = {2 q \over q^2 + 1}.
$$
The expression (5.20) has the same limit as (3.8) when $q \rightarrow 1$, 
as it must for the right-hand side of (5.10) to exist. With $\sigma$
given by (5.7), the covariant derivative is compatible with the metric
(5.13) and torsion free.  Comparing (3.18) with (3.8) one sees that one
obtains for the second example (3.11) a covariant derivative compatible
with the metric (5.13) by the replacement $q \mapsto q^{-4}$ in (5.20).
The dependence on $q$ comes through the conditions (5.19). Since 
$S(q) = - S(-q^{-1})$ the same matrix serves two different values of the
parameter $q$.
 
The metric we have chosen is not symmetric with respect to $\sigma$.
That is
$$
g^{ab} \neq S^{ab}{}_{cd} g^{cd}
$$
in general. If one wishes to find a metric symmetric in the above sense
then one must consider (5.16) as an equation for $S$ and the metric and
add the additional equation
$$
g^{ab} = S^{ab}{}_{cd} g^{cd}.                                    \eqno(5.21)
$$
The system (5.16), (5.21), without the restricton we have placed on the
coefficients $g^{ab}$, if it has a solution, would yield a symmetric
metric with a compatible connection.

Restricting one's attention to hermitian $x$ and $y$ and real $q$,
in the limit $q \rightarrow 1$ one obtains on the ordinary 2-plane
a metric whose Gaussian curvature $K$ is given by
$$
K_1 = x^2 + y^2,                                     \qquad 
K_2 = x^{-4} (1 + y^4)                                            \eqno(5.22)
$$
respectively for the two examples of Section~3.  This can be calculated
using the $q \rightarrow 1$ limit of (5.10).  It is easy to characterize
all metrics which can be obtained in this way.  In the limit 
$q \rightarrow 1$ the commutator determines a Poisson bracket on the
2-plane given as usual by
$$
\{f,g\} = \lim_{q\rightarrow 1} {1\over q-1} [f,g].               \eqno(5.23)
$$  
Define 
$$
p_a = \lim_{q\rightarrow 1} (q-1) \lambda_a.
$$
In the limit the differential can be written then in the form
$$
df = \{p_a, f\} \theta^a.                                          \eqno(5.24)
$$ 
If we write $\theta^a = \theta^a_b dx^b$ in terms of $dx^a$ from this
it follows that the equation 
$$
\{p_c, x^a\} \theta^c_b = \delta^a_b.                              \eqno(5.25)
$$
must have a solution for $p_a$ polynomial in the variables.  This is not
always the case.  That is, not all metrics with polynomial curvature can
be obtained as were those given by (5.22). For example consider the flat
metric $\theta^a_b = \delta^a_b$. The equations (5.25) become the
equations $\{p_a, x^b\} = \delta^a_b$.  Using the expression for the
Poisson bracket for the generalized quantum plane, $\{x,y\} = xy$, one
sees immediately that there is no solution for the $p_a$.

The generalized quantum plane has two outer derivations defined by
$$
e^{(0)}_1 x = x, \qquad e^{(0)}_1 y = 0, \qquad 
e^{(0)}_2 x = 0, \qquad e^{(0)}_2 y = y.                         \eqno(5.26)
$$
The corresponding $\theta^a$ are given by
$$
\theta^1 = x^{-1} dx, \qquad \theta^2 = y^{-1} dy.               \eqno(5.27)
$$
Our construction yields then the ordinary flat metric. If one were
to extend the algebra to the Heisenberg algebra then these derivations
would become interior. To obtain a metric which is almost flat one can
add to (5.26) a `small' inner derivation of the form given in Section~2
but using $\lambda_a$ which are `small' of the order of some expansion
parameter $\epsilon$. One defines
$$
e_a = e^{(0)}_a + \epsilon \, {\rm ad}\, \lambda_a               \eqno(5.28)
$$
and proceeds as above but retaining only corrections of first order in
$\epsilon$.  A problem closely related to this has been examined in
another context by one of the authors (Madore 1995b).

The equations (5.17), (5.19) admit also the solution
$$ 
S = {1 \over q^2 + 1}\pmatrix{
- 2 q   & 0       & 0       & 1 - q^2   \cr
0       & 1 - q^2 & 2 q     & 0         \cr
0       & 2 q     & q^2 - 1 & 0         \cr
q^2 - 1 & 0       & 0       & - 2 q
}                                                                 \eqno(5.29)
$$
but the corresponding covariant derivative defined by (5.10) is singular
in the limit $q\rightarrow 1$.

A complete study of the matrix case has not been made. However for the
particular case (5.12) it is easy to see that the unique torsion-free 
covariant derivative compatible with the metric (5.13) is given by
$$
D\theta^a = - {1\over 2} C^a{}_{bc} \theta^b \otimes \theta^c.   \eqno(5.30)
$$
The ordinary quantum plane with the differential calculus given by
(3.13) has no metric connection but it posesses a unique 1-parameter
family of linear connections whose curvature is polynomial in the
variables $x$ and $y$ (Dubois-Violette {\it et al.} 1995a).
The precise property of the curvature $K_2$ in (5.22) which associates
the corresponding metric to the $GL_q(2)$-invariant differential
calculus (3.13) is not clear.  We refer to Madore \& Mourad (1996) for a
description of the possible relevance to the theory of gravity of the
relation between linear connections on the one hand and commutation
relations on the other.

\section{Conclusions}
We have shown that each differential calculus and set of commutation
relations between the 1-forms and the elements of the algebra gives rise
in the case of the generalized quantum plane to a metric connection
which remains regular in the limit $q \rightarrow 1$. Not all metrics
with polynomial curvature can be obtained in this way.

\section{Acknowledgment} One of the authors (J.M.) would like to thank M.
Dubois-Violette, T. Masson and J. Mourad for interesting conversations.
%\vskip 1cm
%\vfill\eject
\parskip 7pt plus 1pt
\parindent=0cm

\section*{References}

Baehr H.C., Dimakis A., M\"uller-Hoissen F. 1995, {\it Differential
calculi on commutative algebras}, J. Phys. A {\bf 28} 3197.

Bresser K., M\"uller-Hoissen F., Dimakis A., Sitarz A. 1995, 
{\it Noncommutative Geometry of Finite Groups}, Preprint GOET-TP 95/95.

Connes A. 1986, {\it Non-Commutative Differential Geometry}, Publications
of the Inst. des Hautes Etudes Scientifique. {\bf 62} 257.

--- 1994, {\it Noncommutative Geometry}, Academic Press.

Cuntz A., Quillen D. 1995, {\it Algebra extensions and nonsingularity},
J. Amer. Math. Soc. {\bf 8} 251.

Dimakis A., M\"uller-Hoissen F. 1992, {\it Noncommutative differential 
calculus, gauge theory and gravitation}, Preprint GOE-TP 33/92.

Dimakis A., Tzanakis C. 1995, {\it Noncommutative Geometry and Kinetic
Theory of Open systems}, hep-th/9508035, J. Phys. A. (to appear).

Dubois-Violette M. 1988, {\it D\'erivations et calcul diff\'erentiel 
non-com\-mutatif}, C. R. Acad. Sci. Paris {\bf 307} S\'erie I 403.

Dubois-Violette M., Kerner R., Madore J. 1989, {\it Gauge bosons in a
noncommutative geometry}, Phys. Lett. {\bf B217} 485; {\it Classical
bosons in a noncommutative geometry}, Class. Quant. Grav. {\bf 6} 1709.

Dubois-Violette M., Madore J., Masson T., Mourad J. 1995a, {\it Linear
Connections on the Quantum Plane}, Lett. Math. Phys. {\bf 35} 351.

--- 1995b, {\it On Curvature in Noncommutative Geometry}, 
Preprint LPTHE Orsay 95/63; q-alg/9512004.

Dubois-Violette M., Michor P. 1995, {\it Connections on Central
Bimodules}, Preprint LPTHE Orsay 94/100.

Dubois-Violette M., Masson T. 1995, {\it On the First-Order Operators
in Bimodules}, Preprint LPTHE Orsay 95/56.

Karoubi M. 1981, {\it Connections, courbures et classes
caract\'eristiques en $K$-th\'eorie alg\'ebrique.}, Current trends in
algebraic topology, Part I, London, Ont.

Koszul J.L. 1960, {\it Lectures on Fibre Bundles and Differential Geometry},
Tata Institute of Fundamental Research, Bombay.

Madore J. 1995a, {\it An Introduction to Noncommutative Differential Geometry
and its Physical Applications}, Cambridge University Press.

--- 1995b, {\it Linear Connections on Fuzzy Manifolds}, 
Preprint LPTHE Orsay 95/42, ESI Vienna 235, hep-th/9506183.

Madore J., Mourad. J. 1996, {\it Noncommutative Kaluza-Klein Theory},
Lecture given at the $5^{\rm th}$ Hellenic School and Workshops on
Elementary Particle Physics (to be published).

Madore J., Masson T., Mourad J. 1995, {\it Linear Connections on 
Matrix Geometries}, Class. Quant. Grav. {\bf 12} 1429.

Mourad. J. 1995, {\it Linear Connections in Non-Commutative Geometry},
Class. Quant. Grav. {\bf 12} 965.

Wess J., Zumino B. 1990, {\it Covariant Differential Calculus on the
Quantum Hyperplane} Nucl. Phys. B (Proc. Suppl.) {\bf 18B} 302.

\end{document}